\documentclass[12pt,dvips]{amsart}
\usepackage{euler, amsfonts, amssymb, latexsym, epsfig,epic}

\newcommand\iso{{\,\cong\,}}

\newtheorem{Theorem}{Theorem} 
\newtheorem{Proposition}{Proposition} 
\newtheorem{Lemma}{Lemma}

\newtheorem*{Corollary}{Corollary}
 
\newtheorem*{Theorem*}{Theorem}
\newtheorem*{Corollary*}{Corollary}
\theoremstyle{remark}

\newcommand\into{\operatorname*{\hookrightarrow}}
\newcommand\Ex{{\em Example. }}

\newcommand\reals{{\mathbb R}}
\newcommand\complexes{{\mathbb C}}
\newcommand\integers{{\mathbb Z}}

\newcommand\naturals{{\mathbb N}}

\theoremstyle{plain}

\newcommand\dfn{\bf} 

\newcommand\GLn{{{GL_n(\complexes)}}}

\newcommand\covers{\succ}

\begin{document}
\pagestyle{plain}

\title{A Schubert calculus recurrence from\\
 the noncomplex $W$-action on $G/B$}
\author{Allen Knutson}
\thanks{AK was supported by an NSF grant and the Sloan Foundation.}
\email{allenk@math.berkeley.edu}
\date{\today}

\begin{abstract}
For $K$ a compact Lie group (with a chosen maximal torus $T$) 
and $G$ its complexification (with a chosen Borel subgroup $B$), the
diffeomorphism $K/T \iso G/B$ lets one see a noncomplex {\em right}
action of the Weyl group on this complex manifold. 

We calculate the action of simple reflections from $W$ on the
cohomology ring, in the basis of Schubert classes, and use it to
give a (nonpositive) recurrence on the structure constants. 

Our main computational tool is equivariant cohomology, which lets one
model cohomology classes by lists of polynomials \cite{A,KK}.
\end{abstract}

\maketitle

\section{Statement of results}\label{sec:alg}

Our main result is a recurrence relation on the structure constants
for multiplication in the (torus-equivariant) cohomology ring of a 
generalized flag manifold $G/B$.
These structure constants, known as ``equivariant Schubert calculus'',
are known to be positive (for the equivariant statement, see \cite{G}). 
Using this recurrence and the ``descent-cycling'' results from \cite{K}, 
we give in this section a recursive algorithm to compute the
structure constants (for $G$ finite-dimensional).

Our recurrence is not manifestly positive, alas, though in practice
it frequently has no negative terms (see the comment on ``anti-Grassmannian
permutations'' below).
Moreover, unlike computations based on representatives for the
Schubert classes (e.g. Schubert polynomials, or as in \cite{B}), this
algorithm does not require us to compute a full product of two
Schubert classes in order to extract just one term.

We briefly recall from \cite{KK} the standard notation we need
(with more detail in the next section): the equivariant
Schubert classes $\{S_w\}$ are indexed by the Weyl group $W$,
whose (strong) Bruhat order is denoted $\geq$ with
covering relation $\covers$. Write
$$ S_w S_v = \sum_u c_{wv}^u S_u, \quad c_{wv}^u \in H^*_T(pt).$$
This ring $H^*_T(pt)$ is 
the polynomial ring in the simple roots of $G$ (for $G$ adjoint),
so any root is naturally an element (in $H^2_T(pt)$ by convention).
The case $l(u) = l(w)+l(v)$ (where $l$ denotes the length function on
the Coxeter group $W$) is called the {\dfn ordinary (cohomology) case},
as those are the only nonvanishing $c_{wv}^u$ in ordinary, non-equivariant,
Schubert calculus.

We now give a recursive algorithm for computing any one $c_{wv}^u$, 
$w$ not the longest element $w_0$, 
in terms of $c_{w'v'}^{u'}$ with $w'>w$ and/or $v'<v$ 
(in the ordinary case, the latter is not necessary).
Let $r=r_\alpha$ be a reflection through a simple root $\alpha$
such that $wr>w$; in the case $W=S_n$ and
$r = (i\ i+1)$ this corresponds to $w$ having an ascent between the
$i$th and $(i+1)$st positions. (There will exist such an $r$ unless 
$w = w_0$.) 
Now consider the two conditions $vr>v$ and $ur<u$.
  \begin{enumerate}
  \item If both $vr>v$ and $ur<u$ hold, then $c_{wv}^u = 0$.
    This simple result is called ``dc-triviality'' in \cite{K}.
  \item If exactly one is true, we can replace the Schubert problem by
    another with a higher $w$:
    \begin{itemize}
    \item If $ur<u$ but $vr<v$, then $c_{wv}^u = c_{wr,vr}^{u}.$
    \item If $vr>v$ but $ur>u$, then $c_{wv}^u = c_{wr,v}^{ur}.$
    \end{itemize}
    This symmetry was introduced in \cite{K} and is called ``descent-cycling''.
  \item If neither is true -- so $vr<v$ and $ur>u$ --
    then we use the new recurrence (theorem \ref{thm:recurrence}):
$$ c_{wv}^u 
 = c_{wr,v}^{ur} + c_{wr,vr}^u - (w\cdot \alpha) c_{w,vr}^u 
  + \sum_{w'\covers w,\, w'\neq wr,\, w' = wr_\beta}
  \langle \alpha,\beta \rangle\, c_{w',vr}^u. $$
  (Here $\beta$ is a positive root, not necessarily simple,
  and $r_\beta$ the reflection through it.)

  The first two terms on the right-hand side are more Schubert structure
  constants. The third is a structure constant times an element of
  $H^2_T(pt)$ (that is negative in the sense of \cite{G}); this term
  drops out when computing ordinary Schubert calculus.
  The coefficients $\langle \alpha,\beta \rangle
  := (\alpha - r_\beta \alpha)/{\beta} \in\integers$ in the last term
  are the other possible source of nonpositivity in this recurrence.
  \end{enumerate}
The algorithm above applies whenever $w\neq w_0$, and writes a Schubert
structure constant in terms of structure constants with higher $w$ and/or 
lower $v$. Therefore it terminates (assuming that $G$ is finite-dimensional),
and only requires that we be able to compute the $\{c_{w_0,v}^u\}$.
These vanish unless $u=w_0$. In the ordinary case, 
the only one met is $c_{w_0,1}^{w_0} = 1$.

Our examples will be in $W=S_n$, where we can speak in terms of
ascents and descents, and the coefficients $\langle \alpha,\beta \rangle$
are always $\pm 1$, as follows. If $\alpha = x_i - x_{i+1}$, 
each $w' = wr_\beta$ with nonvanishing $\langle \alpha,\beta \rangle$
agrees with $w$ except in two places, exactly one of which
is the $i$th or $(i+1)$st place. If the two places switched straddle
the $i,i+1$ divide, the coefficient $\langle \alpha,\beta \rangle$ 
is $+1$; if they are both to one side of it, the coefficient is $-1$.

\Ex Let $G=GL_4(\complexes)$, $W=S_4$, and $w=1234, v=u=2413$ (so of course
the answer is $c_{wv}^u = 1$, an ordinary case). 
We first apply the algorithm with $r=(23)$
(for illustrative purposes only; in all other cases to follow we will use
the $r=(i,i+1)$ with {\em least} $i$ such that $w_i < w_{i+1}$).
It tells us to descent-cycle:
$$ c_{12|34,24|13}^{24|13} = c_{1324,2143}^{2413} $$
In this and all subsequent examples, we put $|$s in the permutations
to indicate the next choice of $r$ used.

Now it is impossible to cycle another descent directly into $w$, so we
apply the recurrence with $r=(12)$:
$$ c_{1|324,2|143}^{2|413} 
 = c_{3124,2143}^{4213} 
 + c_{3124,1243}^{2413}
 - c_{1423,1243}^{2413}
 + c_{2314,1243}^{2413}
$$
(The equivariant term would be $(y_1-y_3) c_{1324,1243}^{2413}$,
but since this is an ordinary-case calculation, that term vanishes for
degree reasons.) The first three of these terms die:
$$ c_{31|24,21|43}^{42|13} = 0, \qquad\qquad
   c_{31|24,12|43}^{24|13} = 0, \qquad\qquad
  -c_{1|423,1|243}^{2|413} = -c_{41|23,12|43}^{42|13} = 0 $$
The fourth term requires the recurrence again, after a descent-cycling:
$$ c_{2|314,1|243}^{2|413} =
  c_{321|4,124|3}^{421|3} 
= c_{3241,1243}^{4231}
+ c_{3241,1234}^{4213} 
+ c_{3412,1234}^{4213} 
+ c_{4213,1234}^{4213} 
$$
These terms each simplify quickly, giving us the desired answer $1$:
$$ c_{32|41,12|43}^{42|31} = c_{3|421,1|243}^{4|321} = 0, \qquad \qquad
   c_{32|41,12|34}^{42|13} = 0, \qquad \qquad
   c_{3|412,1|234}^{4|213} = 0 $$
$$ c_{421|3,123|4}^{421|3} 
 = c_{42|31,12|34}^{42|31} = c_{4321,1234}^{4321} = 1. $$

Call a permutation $w\in S_n$ {\dfn anti-Grassmannian} if $w$ has at most one
ascent, i.e. if $w_0 w$ is a Grassmannian permutation. Note that if $w$
is anti-Grassmannian, and we're in the ordinary case, then this
recurrence relation has no negative terms.

\Ex Let $w=532164, v=132546, u=645231$. This ordinary $c_{wv}^u$ is
$2$ (this and all other experiments were done with \cite{ACE}). We
cannot descent-cycle into $w$, so we apply the recurrence:
\begin{align*}
  & c_{5321|64,1325|46}^{6421|53} \\
&= c_{532|614,132|546}^{642|513}
&+& c_{532|614,132|456}^{642|153}
&+& c_{6321|54,1324|56}^{6421|53}
&+& c_{5|62134,1|32456}^{6|42153}
&+& c_{53|6124,13|2456}^{64|2153}
&+& c_{532|461,132|456}^{642|153} \\
&= c_{536241,132546}^{645231} 
&+& 0
&+& c_{63251|4,13245|6}^{64251|3} 
&+& 0
&+& c_{56312|4,12345|6}^{64215|3} 
&+& 0 \\
&= c_{536241,132546}^{645231} 
&+& 0
&+& c_{632|541,132|465}^{642|531} 
&+& 0
&+& 0
&+& 0 \\
\end{align*}
Then
\begin{align*}
   c_{53|6241,13|2546}^{64|5231} 
&=& c_{5|63241,1|32546}^{6|54231} 
&+& c_{5|63241,1|23546}^{6|45231} 
&+& c_{63|5241,12|3546}^{64|5231} 
&+& c_{54|6231,12|3546}^{64|5231} \\
&=& 0
&+& 0
&+& c_{653241,123546}^{654231} 
&+& c_{5|64231,1|23546}^{6|54231} \\
&=& 0
&+& 0
&+& c_{653241,123546}^{654231} 
&+& 0
\end{align*}
\begin{align*}
   c_{6532|41,1235|46}^{6542|31} 
&=& c_{653|421,123|546}^{654|321} 
&+& c_{653|421,123|456}^{654|231}
&+& c_{6542|31,1234|56}^{6542|31} \\
&=& 0
&+& 0
&+& c_{654321,123456}^{654321} &=& 1
\end{align*}
and
\begin{align*}
 c_{6352|41,1324|56}^{6452|31} 
&=& c_{63|5421,13|2456}^{64|5321} \\
&=& c_{653|421,132|456}^{654|321} 
&+& c_{653|421,123|456}^{645|321} 
&+& c_{64|5321,12|3456}^{64|5321} \\
&=& 0
&+& 0
&+& c_{654321,123456}^{654321} &=&1.
\end{align*}
In all, $2=1+1$. Each time we used the recurrence in this example, $w$ was
anti-Grassmannian, and so there were no minus signs.

In equivariant cohomology the base case $c_{w_0,v}^{w_0}$ is harder, 
but the results of \cite{B,W,GK} can be adapted to compute it 
(theorem \ref{thm:Sara} in the next section).

\Ex Let $W=S_3$, $w=231, v=213, u=231$. This $c_{231,213}^{231} = y_2-y_1$, and
the algorithm computes it (non-\cite{G}-positively) as
\begin{align*}
 c_{2|31,2|13}^{2|31} 
&=& c_{321,213}^{321}
&+& c_{321,123}^{231} 
&-& (y_3 - y_2) c_{2|31,1|23}^{2|31} \\
&=& (y_3 - y_1) 
&+& 0
&-& (y_3 - y_2) c_{321,123}^{321}.
\end{align*}

\section{The Schubert basis of $T$-equivariant cohomology of $G/B$}

We set up our conventions, and include some standard material on
equivariant Schubert calculus from \cite{A,KK}.

\newcommand\codim{{\rm codim}}
Fix a pinning $(G,B,B_-,T^\complexes,W)$ of a complex Lie group;
our motivating example is $G=\GLn$, $B$ the upper triangulars, $B_-$
the lower triangulars, $T^\complexes$ the diagonals, and $W \iso S_n$.
For each element $w\in W$, the {\dfn Schubert cycle} 
$X_w$ is the orbit closure $\overline{B_- wB}$. Being a $T$-invariant cycle,
it induces an element of $T$-equivariant cohomology; the degree of 
this cohomology class is $\codim_\reals X_w = 2\, l(w)$, twice the length
of the Weyl group element $w$. The forgetful map
$H^*_T(G/B) \to H^*(G/B)$ takes this equivariant Schubert class to the
ordinary one, so it is no harm to work in this richer cohomology theory.
The equivariant Schubert classes are again a basis for cohomology,
but over the base ring $H^*_T(pt)$, which (for $G$ adjoint) is just
the polynomial ring in the simple roots (each formally given degree $2$).

The pullback ring homomorphism 
from $H^*_T(G/B) \to H^*_T((G/B)^T) = \bigoplus_W H^*_T(pt)$
takes an equivariant class $\Psi$ to $\{ \Psi|_w \}_{w\in W}$,
an $H^*_T(pt)$-valued function on $W$. If $\Psi = S_w$, then the
support of this function is $\{v\in W: v\geq w\}$. This upper
triangularity implies that the pullback map is $1:1$. Accordingly,
we will do all our calculations with these lists of polynomials.

In \cite{A} was calculated the image of this restriction map:
a list $\{\alpha|_v\}$ comes from a cohomology class if and only if
$$ \forall v\in v, \beta\in\Delta, \quad \alpha|_v - \alpha|_{r_\beta v}
\hbox{ is a multiple of }\beta. $$
These conditions are nowadays viewed in the more general framework 
of \cite{GKM}.
Hereafter we define a {\dfn class} $p$ to be a list $\{p|_w\}_{w\in W}$
of elements of $H^*_T(pt)$, satisfying these {\dfn GKM conditions}.

We also recall a characterization of the Schubert class $S_w$:
it is homogeneous of degree $2l(w)$, its
restriction $S_w|_v$ vanishes at $v$ shorter or of the same length as $w$
(except at $w$ itself),
and its restriction $S_w|_w$ at its bottom point is 
$$ S_w|_w = \prod_{\beta\in \Delta_+,\, r_\beta w<w} \beta. $$
Plainly the GKM conditions force a class vanishing below $w$ to be
a multiple of this monomial at $w$; the Schubert class is characterized by
this multiple being $1$.

At this point we abandon the geometry and work only with this
combinatorial model of the equivariant cohomology ring, and its
basis of Schubert clases, much as in \cite{KT}. Accordingly, while we
will relate our constructions to geometry wherever possible, we will
only give proofs of the combinatorial statements.

\begin{Lemma}
  The restriction $S_v|_w$ of the class $S_v$ to a point $w$ is an
  equivariant Schubert structure constant, $c_{wv}^w$.
\end{Lemma}

\begin{proof}
  The class $S_w S_v$ vanishes when restricted to $u$ unless $u\geq w$
  (since already $S_w$ does), 
  so the upper triangularity tells us that $c_{wv}^w = (S_w S_v)|_w / S_w|_w$.
\end{proof}

\begin{Theorem}\cite{B,W,GK}\label{thm:Sara}
  Let $I$ be a reduced expression for $w\in W$, whose $i$th entry is the
  reflection $r_i$ through the simple root $\alpha_i$. Then for each $v\in W$,
$$ S_v|_w = \sum_{J\subseteq I} 
        \prod_I \big( \hat\alpha_i^{[i\in J]} r_i \big) \cdot 1$$
  where the sum is taken over reduced subwords $J$ with product $v$, 
  and the $\hat\alpha_i$ are multiplication operators only included in the
  ordered product if $i\in J$.
\end{Theorem}

\Ex If $I = r_{12} r_{23} r_{12}$ is a reduced word for $321\in S_3$,
then 
$$ S_{213}|_{321} 
 = \widehat{y_2-y_1} r_{12} r_{23} r_{12}\cdot 1 
 + r_{12} r_{23} \widehat{y_2-y_1} r_{12}\cdot 1 
 = \widehat{y_2-y_1} \cdot 1 
 + r_{12} r_{23} \cdot (y_2-y_1) = (y_2-y_1) + (y_3-y_2). $$
Whereas if we use $I = r_{23} r_{12} r_{23}$, we'd get
$$ S_{213}|_{321} 
 = r_{23} \widehat{y_2-y_1} r_{12} r_{23} \cdot 1 
 = r_{23} \cdot (y_2-y_1) = y_3 - y_1. $$

Combining this lemma and theorem, 
we have a (\cite{G}-positive) formula for the base case $c_{w_0,v}^{w_0}$.
One special case is worthy of note: in the $W=S_n$ case, using the 
lexicographically last reduced expression for $w_0$, the terms $J$ 
correspond to the {\em rc-graphs} for $v$, and $c_{w_0,v}^{w_0}$ is equal
to the double Schubert polynomial for $v$, evaluated with $x_i$ set equal
to $y_{n+1-i}$.

\section{Left and right divided difference operators}

We first define left and right actions of $W$ on the ring of classes:
$$ (w \cdot p)|_v := w \cdot (p|_{w v}) $$
which uses the action of $W$ on the base ring $H^*_T(pt)$, and
$$ (p \cdot w)|_v := p|_{v w} $$
which is $H^*_T(pt)$-linear.

\begin{Proposition}
  If $p$ is a class and $w\in W$, then $w\cdot p$ and $p\cdot w$
  are classes. Both actions define ring automorphisms, but only the
  second is an $H^*_T(pt)$-algebra automorphism.
\end{Proposition}

\begin{proof}
  We need to check the GKM conditions. For the first,
$$ (w\cdot p)|_v - (w\cdot p)|_{r_\beta v} 
 = w\cdot (p|_{wv} - p|_{w r_\beta v})
 = w\cdot (p|_{wv} - p|_{ r_{w^{-1}\cdot \beta} wv})
$$
  is indeed a multiple of $\beta$, 
  since $p|_{wv} - p|_{ r_{w^{-1}\cdot \beta} wv}$ is a multiple of 
  $w^{-1}\cdot \beta$.
  
  For the second, 
$$ (p\cdot w)|_v - (p\cdot w)|_{r_\beta v} = p|_{vw} - p|_{r_\beta vw} $$
  is even more obviously a multiple of $\beta$.

  The ring automorphism statement is obvious. The first fails to be an
  algebra automorphism, exactly because $W$ is acting on the coefficients.
\end{proof}

There is a geometric reason for this proposition; the first action arises from 
the left action of $N(T)$ on $G/B$, by diffeomorphisms that {\em normalize} 
the $T$-action, and thereby induce ring automorphisms of $T$-equivariant 
cohomology, whereas the second action arises from the right action 
of $W$ on $G/B \iso K/T$,
which {\em commutes} with the left $T$-action (and even $K$-action),
and thereby induce algebra automorphisms.

Note that the right action does not preserve the complex structure on
$G/B$, and does not in general take a Schubert class to a positive
combination of other Schubert classes.

For each simple root $\alpha$, 
define the {\dfn left divided difference operator} $\partial_\alpha$ by
$$ \partial_\alpha p := \frac{1}{\alpha}(p - r_\alpha\cdot p). $$
This is a famous degree $-2$ endomorphism of $H^*_T(G/B)$ (though not a
module homomorphism), and variants of it have long since been used 
to inductively construct the Schubert classes, as in parts (3) and (4)
of proposition \ref{prop:ddops} below. For us, it only serves to
motivate the definition below of the right divided difference operators.

For each simple root $\alpha$, 
define the {\dfn Chern class $c_{-\alpha}$}
associated to $\alpha$ by 
$$ c_{-\alpha}|_w := w\cdot (-\alpha). $$
By construction, it satisfies $v \cdot c_{-\alpha} = c_{-\alpha}$ 
for all $v\in W$.
This list $c_{-\alpha}$ is easily seen to satisfy the GKM conditions, 
so is a class; geometrically,
it arises as the equivariant first Chern class of the Borel-Weil line bundle
associated to $-\alpha$.\footnote{%
In the Borel-Weil theory, one generally associates line bundles to
{\em dominant} weights; the strongly dominant weights give ample line bundles.
The negative simple roots are far from being dominant,
and these line bundles are not typically ample.} 
That being a $G$-equivariant line bundle, 
it is $N(T)$-equivariant, explaining $c_{-\alpha}$'s $W$-invariance. 

Note that $c_{-\alpha}$ is not a zero divisor in the ring of classes.
With it, we can define the {\dfn right divided difference
operator} $\partial^\alpha$ by
$$ \partial^\alpha p := \frac{1}{c_{-\alpha}}(p - p\cdot r_\alpha), $$
and this {\em is} an $H^*_T(pt)$-module homomorphism.
Geometrically, this arises (as in \cite{BGG}) from composing pushforward
and pullback of the $T$-equivariant morphism $G/B \to G/P_{\alpha}$.

\begin{Proposition}\label{prop:ddops}
  \begin{enumerate}
  \item The left and right divided difference operators take classes
    to classes.
  \item Left divided difference operators commute with 
    right divided difference operators.
  \item $\partial^{\alpha} S_w = S_{w r_\alpha}$ if $w r_\alpha < w$, $0$
    otherwise.
  \item $\partial_{\alpha} S_w = S_{r_\alpha w}$ if $r_\alpha w < w$, $0$
    otherwise.
  \item $\partial^{\alpha} (pq) 
    = (p - c_{-\alpha} \partial^\alpha p) (\partial^\alpha q) 
    + (\partial^\alpha p) q. $
  \end{enumerate}
\end{Proposition}

Readers wondering about the unfortunate minus sign in $c_{-\alpha}$ can now 
trace it to our desire to have $\partial^{\alpha} S_w = +S_{w r_\alpha}$.

\begin{proof}
  We first show that $\partial_\alpha,\partial^\alpha$ take classes to lists of
  polynomials, not rational functions:
$$ (\partial_\alpha p)|_v = \frac{1}{\alpha}(p|_v - p|_{r_\alpha v}), \qquad
 (\partial^\alpha p)|_v = \frac{1}{c_{-\alpha}|_v}(p|_v - p|_{v r_\alpha})
 = \frac{1}{-v \cdot \alpha}(p|_v - p|_{r_{v\cdot\alpha} v}) $$
  In both cases, these are polynomials because $p$ satisfies the 
  GKM conditions.

  Next we look at $(\partial_\alpha p)|_v - (\partial_\alpha p)|_{r_\beta v}$
$$ = \frac{1}{\alpha}
  (p|_v - r_\alpha \cdot (p|_{r_\alpha v}) - 
p|_{r_\beta v} + r_\alpha \cdot (p|_{r_\alpha r_\beta v})) 
=
 \frac{1}{\alpha} (
  (p|_v - p|_{r_\beta v}) -
r_\alpha \cdot (p|_{r_\alpha v} - p|_{r_\alpha r_\beta v}))
$$
  If $\alpha = \beta$, then this is zero. Otherwise $gcd(\alpha,\beta)=1$,
  so we can ignore the $1/\alpha$ when testing divisibility by $\beta$,
  and each the terms inside is a multiple of $\beta$. So $\partial_\alpha p$
  satisfies the GKM conditions.

  Commutativity follows from the fact that $r_\alpha \cdot c_\beta = c_\beta$
  and the fact that left and right multiplication in $W$ commute.

  We can use commutativity to prove that $\partial^\alpha p$ is a class:
$$ (\partial^\alpha p)|_v - (\partial^\alpha p)|_{r_\beta v} 
= \beta (\partial_\beta \partial^\alpha p)|_v
= \beta (\partial^\alpha \partial_\beta p)|_v
$$
  We already know $\partial_\beta p$ is a class, and $\partial^\alpha$ of
  a class is a list of polynomials, therefore this difference is a
  multiple of $\beta$.

  Now consider $\partial^\alpha S_w$. This vanishes except at 
  $\{v : v\geq w \hbox{ or } v r_\alpha \geq w \}$.
  By the upper triangularity of the
  Schubert classes, this class is an $H^*_T(pt)$-linear combination 
  of $\{ S_v : v\geq w \hbox{ or } v r_\alpha \geq w \}$.
  All of these $S_v$ have degree
  higher than that of $\partial_\alpha S_w$, unless $v = w r_\alpha < w$.
  We check the value at $w r_\alpha$:
$$ (\partial^\alpha S_w)|_{wr_\alpha} 
= \frac{1}{c_{-\alpha}} (S_w - S_w\cdot r_\alpha)|_{wr_\alpha}
= \frac{1}{w\cdot \alpha} (S_w|_{wr_\alpha} - S_w|_w) 
= \frac{S_w|_w}{-w\cdot\alpha} = S_{wr_\alpha}|_{wr_\alpha} 
$$
  So by the aforementioned characterization, 
  $\partial^\alpha S_w = S_{w r_\alpha}.$
  Whereas if $w r_\alpha > w$, there are no possible terms and 
  $\partial^\alpha S_w = 0$. 
  The proof for $\partial_\alpha$ is exactly the same.
  Alternately, one can prove it for the special case $w=w_0$ and
  use the commutativity.

  Finally, we compute
$$ \partial^\alpha (pq) = \frac{1}{c_{-\alpha}}(pq - (pq)\cdot r_\alpha) 
 = \frac{1}{c_{-\alpha}}(pq - p(q\cdot r_\alpha) + p(q\cdot r_\alpha)
 - (p\cdot r_\alpha) (q\cdot r_\alpha) ) $$
$$ = p (\partial^\alpha q) + (\partial^\alpha p) (q\cdot r_\alpha) 
   = p (\partial^\alpha q) 
   + (\partial^\alpha p) (q - c_{-\alpha} \partial^\alpha q) 
   = (p - c_{-\alpha} \partial^\alpha p)(\partial^\alpha q) 
   + (\partial^\alpha p) q $$
as claimed.
\end{proof}

\subsection{On the left vs. right symmetry.}
We comment on a deeper algebraic reason, itself from a more trivial
geometric reason, for the apparent symmetry between the left and right
actions of $W$. The base ring $H^*_T(pt) \into H^*_T(G/B)$ 
can be characterized as the
right-$W$-invariant subring of $H^*_T(G/B)$, and is generated freely by 
the simple roots. There is another subring $C \leq H^*_T(G/B)$
of {\em left}-$W$-invariant classes, 
generated freely by the Chern classes $\{c_{-\alpha}\}$, 
and in particular there is
a natural isomorphism $\phi: H^*_T(pt) \iso C$ taking 
$p \mapsto c_p = \{c_p|_v := v\cdot p\}$.  The Schubert classes
$\{S_w\}$ are a $C$-basis also for $H^*_T(G/B)$ (proved using the same
upper triangularity), and we can consider the structure constants in
$$ S_w S_v = \sum_u d_{wv}^u S_u, \quad d_{wv}^u \in C.$$
Then the best statement is that $\phi(c_{wv}^u) = d_{w^{-1}v^{-1}}^{u^{-1}}$.

We will not need or derive this:
the most natural derivation is via ``double Schubert calculus'', 
in which ones multiplies the cohomology classes of $G$-orbit closures on 
$(G/B)^2$, and then the symmetry comes from switching the two $G/B$ factors.
We found it simpler to work with a single $G/B$, at the expense of 
making the symmetry more mysterious.

\section{Multiplying by Chern classes}

Perhaps the main new idea in this paper 
is a formula for $c_{-\alpha} S_w$. First a lemma:

\begin{Lemma}\label{lem:cover}
  Let $w' = wr_\beta$ cover $w$ in the strong Bruhat order, 
  $\beta \in \Delta_+$. Then
  $$ S_{w'}|_{w'} = S_w|_{w'} (w\cdot \beta). $$
\end{Lemma}

This is quite straightforward using $\partial^\beta$ if $\beta$ is a
simple root (i.e. if the covering relation is in the {\em weak} Bruhat
order); however we need it for all positive roots.

\begin{proof}
  Fix $I$ a reduced word for $w'$. Since $w' \covers w$, there exists a letter
  $b\in I$ such that $I\setminus\{b\}$ is a reduced word for $w$.
  Multiplying out, we can write $w' = w_1 r_b w_2$, and $w = w_1 w_2$. 

  We prove now that such a $b$ is unique (surely a basic fact about the
  Bruhat order, but one for which we lack a reference).
  If $b$ is not unique, we can break $I$ up into $I_1 b I_2 b' I_3$,
  with $w = I_1 b I_2 I_3 = I_1 I_2 b' I_3$. So $b I_2 = I_2 b'$, and
  $b I_2 b' = I_2$, and finally $I_1 b I_2 b' I_3 = I_1 I_2 I_3$, 
  contradicting its reducedness.

  Now we apply theorem \ref{thm:Sara}: 
$$ S_{w'}|_{w'} = \prod_I \big( \hat\alpha_i r_i) \cdot 1,\qquad
   S_w   |_{w'} = \prod_I \big( \hat\alpha_i^{[i\neq b]} r_i) \cdot 1 $$
Since $w' = w w_2^{-1} r_b w_2 = w r_{w_2^{-1}\cdot \alpha_b}$, we have
$$ \frac{S_{w'}|_{w'}}{S_w|_{w'}} = w_1 \cdot \alpha_b 
        = w w_2^{-1} \cdot \alpha_b = w \cdot (w_2^{-1}\cdot \alpha_b) $$
as was to be shown.
\end{proof}

\begin{Proposition}
  Let $\alpha$ be a simple root, $c_{-\alpha}$ the corresponding Chern class,
  and $S_w$ a Schubert class. Then
$$ c_{-\alpha} S_w 
  = -(w\cdot \alpha) S_w + \sum_{w' \covers w,\, w' = w r_\beta} 
        \langle \alpha,\beta \rangle S_{w'}. $$

\end{Proposition}

\begin{proof}
  By degree considerations and the upper triangularity, 
$$ c_{-\alpha}\, S_w = c_{-\alpha}|_w \, S_w + \sum_{v \covers w} d_{v} S_{v}$$
  for some $\{d_{v} \in \integers\}$. Restricting to $w$, all the $v$-terms
  vanish, and we see the coefficient on $S_w$ is as claimed.
  Restricting now to a point $w' \covers w$,
$$ c_{-\alpha}|_{w'}\, S_w|_{w'} 
  = c_{-\alpha}|_w \, S_w|_{w'} +  d_{w'} S_{w'}|_{w'} $$
Dividing by $S_w|_{w'}$ (and using lemma \ref{lem:cover}), this is 
$$ -w'\cdot\alpha = -w\cdot\alpha + d_{w'} (w\cdot\beta) $$
for $w' = wr_\beta$. We can further rewrite as
$$ d_{w'} = w\cdot \frac{\alpha - r_\beta \alpha}{\beta} $$
which is essentially the definition of $\langle \alpha,\beta \rangle$.
\end{proof}

In particular, as claimed in the Abstract,
we now have a formula for $S_w \cdot r_\alpha$:

\begin{Corollary*}
  If $w r_\alpha > w$, then $S_w \cdot r_\alpha = S_w$; otherwise
$$ S_w \cdot r_\alpha 
        = S_w - (w\cdot \alpha) S_{wr_\alpha} -
        \sum_{w' \covers wr_\alpha, \, w' = w r_\alpha r_\beta} 
        \langle \alpha,\beta \rangle S_{w'}. $$ 
\end{Corollary*}

It is equivalent and in fact more convenient to work with the
$\partial^\alpha$ directly.

\section{The recurrence relation}

In this section we prove the recurrence relation on the
Schubert structure constants $c_{wv}^u \in H^*_T(pt)$ in
$ S_w S_v = \sum_u c_{wv}^u S_u$.

Fix a simple root $\alpha$, and let $r = r_\alpha$ be the simple reflection
through it. We introduce a 
convenient {\dfn underline/overline convention}, where a term in a
sum involving $\overline{wr}$ (resp. $\underline{wr}$) only contributes
if $wr > w$ (resp. $wr < w$). For example, we have the single equation
$$ \partial^\alpha S_w = S_{\underline{wr}} $$
encompassing both cases
$ \partial^\alpha S_w = S_{wr} $ if $wr<w$,
$ \partial^\alpha S_w = 0$ if $wr>w$.

The descent-cycling lemmata from \cite{K} are also concisely
expressed with this convention: if $w<wr, v<vr$, then
$$ c_{wv}^u = c_{wr,v}^{\overline{ur}} = c_{w,vr}^{\overline{ur}}. $$
If $ur>u$, this is descent-cycling; if $ur<u$, it is dc-triviality.

\begin{Theorem}\label{thm:recurrence}
  Let $u,v,w \in W$, and $r=r_\alpha$ a simple reflection
  such that $ur>u, vr>v, wr>w$. Then
$$ c_{w,vr}^u 
 = c_{wr,vr}^{ur} + c_{wr,v}^u - (w\cdot \alpha) c_{wv}^u 
  + \sum_{w'\covers w,\, w'\neq wr,\, w' = wr_\beta}
  \langle \alpha,\beta \rangle\, c_{w',v}^u. $$
\end{Theorem}

\begin{proof}
  We apply $\partial^\alpha$ to both sides of the equation 
  $S_w S_v = \sum_u c_{wv}^u S_u$ and equate terms. We do {\em not}   
  yet assume that $w,v,u$ satisfy $ur>u, vr>v, wr>w$.
  (In fact $w,v$ will not match those in the statement of the theorem.)
  On the left side:
$$ \partial^\alpha (S_w S_v) 
 = (S_w - c_{-\alpha} \partial^\alpha S_w) (\partial^\alpha S_v )
 + (\partial^\alpha S_w) S_v 
 = (S_w - c_{-\alpha} S_{\underline{wr}}) S_{\underline{vr}} 
        + S_{\underline{wr}} S_v $$

$$ = \bigg(S_w + (wr\cdot \alpha) S_{\underline{wr}}
   - \sum_{w' \covers \underline{wr},\, w' = w r r_\beta} 
   \langle \alpha,\beta \rangle S_{w'}\bigg) S_{\underline{vr}} 
        + S_{\underline{wr}} S_v $$
$$ = \sum_u S_u \big(
 c_{w,\underline{vr}}^u + (wr\cdot \alpha) c_{\underline{wr},\underline{vr}}^u 
 - \sum_{w' \covers \underline{wr},\, w' = w r r_\beta} 
   \langle \alpha,\beta \rangle c_{w',\underline{vr}}^u 
   + c_{\underline{wr},v}^u \big)$$
(Remember: the underline convention tosses out this sum unless $wr < w$.)

  On the right:
$$ \partial^\alpha \sum_u c_{wv}^u S_u
= \sum_u c_{wv}^u \partial^\alpha S_u
= \sum_u c_{wv}^u  S_{\underline{ur}}
= \sum_u c_{wv}^{\overline{ur}}  S_u $$
Equating coefficients of $S_u$:
$$ 
  c_{w,\underline{vr}}^u 
+ (wr\cdot \alpha) c_{\underline{wr},\underline{vr}}^u 
- \sum_{w' \covers \underline{wr},\, w' = w r r_\beta} 
   \langle \alpha,\beta \rangle c_{w',\underline{vr}}^u 
   + c_{\underline{wr},v}^u 
= c_{wv}^{\overline{ur}} $$
We now consider several special cases. If $vr>v$, then this reduces to
$$ c_{\underline{wr},v}^u = c_{wv}^{\overline{ur}} $$
which is a particularly eloquent summary of several cases of descent-cycling
and dc-triviality. In other words, the formula is known 
in this case, so we will assume $vr<v$.

Descent-cycling also accounts for the $wr<w$ case (easily shown) and
the $ur<u$ case (which is harder). So we will assume that $wr<w$, $vr<v$, 
and $ur>u$.
{\em We now switch the names $w\leftrightarrow wr$,
$v\leftrightarrow vr$,} to make it consistent across variables $u,v,w$,
and also so that the lengths of the Weyl group elements are reflected
in the notation: $l(wr) = l(w) + 1$, etc. 
(One drawback: the $v$ in the algorithm in section \ref{sec:alg}
is the $vr$ of this statement.)

We no longer have need of the 
over/underline convention, since we already know the relative lengths.
$$ 
  c_{wr,v}^u 
+ (w\cdot \alpha) c_{w,v}^u 
- \sum_{w' \covers w,\, w' = w r_\beta} 
   \langle \alpha,\beta \rangle c_{w',v}^u 
   + c_{w,vr}^u 
= c_{wr,vr}^{ur} $$
To go from here to the claimed formula, we move the first three terms over
to the right side of the equation, and pull out the $2c_{wr,v}^u$ term from
the sum ($1$ of which cancels the first term).
\end{proof}

\section{The ordinary-cohomology case}

If $u,v,w \in W$ have $l(u)+l(v)+l(w) = \dim_\complexes G/B$, we can consider
the integral $\int_{G/B} S_u S_v S_w = c_{wv}^{w_0 u} \in \naturals$.
In particular, one sees a three-fold symmetry of the ordinary
Schubert structure constants not visible in the usual definition;
accordingly we denote these integrals by $\{c_{wvu}\}$.

\begin{Corollary}
  Let $u,v,w\in W$ have $l(u)+l(v)+l(w) = \dim_\complexes G/B$, and
  $ur>u,vr>v,wr>w$. Then
$$ c_{w,vr,ur} 
 = c_{wr,vr,u} + c_{wr,v,ur} 
  + \sum_{w'\covers w,\, w'\neq wr,\, w' = wr_\beta}
  \langle \alpha,\beta \rangle\, c_{w',v,ur} $$
\end{Corollary}

\begin{proof}
  This is off from theorem \ref{thm:recurrence} by having switched 
  $u\leftrightarrow ur$ (since in bringing it from a superscript to a
  a subscript, the condition $ur>u$ gets flipped). And the equivariant
  term $(w\cdot \alpha) c_{w,v,ur}$ automatically vanishes, since
  $l(w)+l(v)+l(ur) < \dim_\complexes G/B$. 
 
\end{proof}

\bibliographystyle{alpha}    

\end{document}